\pdfoutput=1
\RequirePackage{ifpdf}
\ifpdf 
\documentclass[pdftex]{sigma}
\else
\documentclass{sigma}
\fi

\begin{document}

\allowdisplaybreaks

\renewcommand{\thefootnote}{$\star$}

\renewcommand{\PaperNumber}{078}

\FirstPageHeading

\ShortArticleName{Frobenius 3-Folds via Singular Flat 3-Webs}

\ArticleName{Frobenius 3-Folds via Singular Flat 3-Webs\footnote{This
paper is a contribution to the Special Issue ``Geometrical Methods in Mathematical Physics''. The full collection is available at \href{http://www.emis.de/journals/SIGMA/GMMP2012.html}{http://www.emis.de/journals/SIGMA/GMMP2012.html}}}

\Author{Sergey I.~AGAFONOV}

\AuthorNameForHeading{S.I.~Agafonov}

\Address{Departmento de Matem\'atica,
Universidade Federal da Paraiba, Jo\~ao Pessoa, Brazil}
\Email{\href{mailto:sergey.agafonov@gmail.com}{sergey.agafonov@gmail.com}}

\ArticleDates{Received May 28, 2012, in f\/inal form October 17, 2012; Published online October 21, 2012}

\Abstract{We give a geometric interpretation of weighted homogeneous solutions to the associativity equation in terms of the web theory and
construct a massive Frobenius $3$-fold germ via a singular $3$-web germ satisfying the following conditions:
1)~the web germ admits at least one inf\/initesimal symmetry,
2)~the Chern connection form is holomorphic,
3)~the curvature form vanishes identically.}


\Classification{53A60; 53D45; 34M35}

\renewcommand{\thefootnote}{\arabic{footnote}}
\setcounter{footnote}{0}

\section{Introduction}

The theory of Frobenius manifolds, having its origin in theoretical physics, has deep interrelations with apparently very dif\/ferent areas of mathematics: Gromov--Witten invariants and quantum cohomology, integrable systems, singularity, deformation of f\/lat connections etc.
We discuss a~new aspect of this fruitful and fast developing theory: its relations with the classical chapter of dif\/ferential geometry, namely the web theory.

\begin{figure}[th]\centering
\includegraphics[width=40mm]{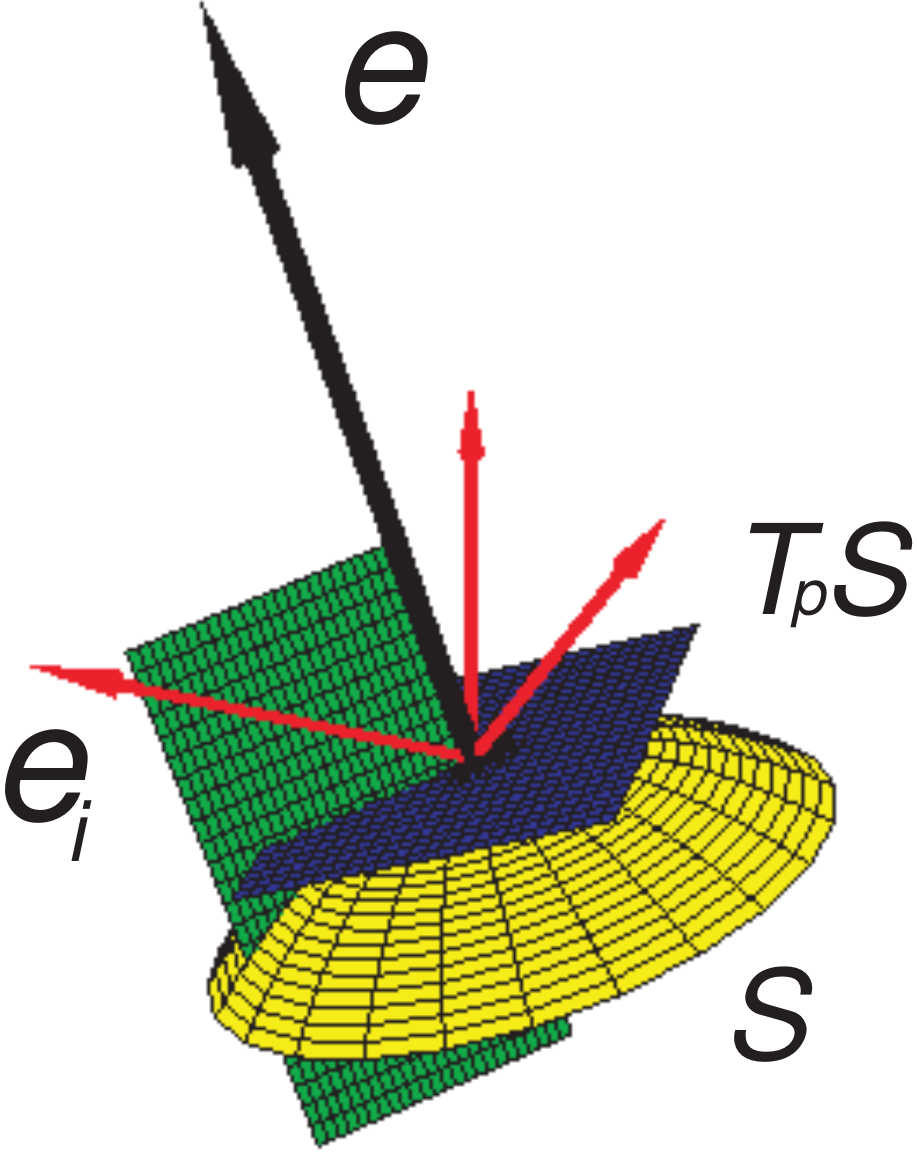} \qquad\quad
\includegraphics[width=40mm]{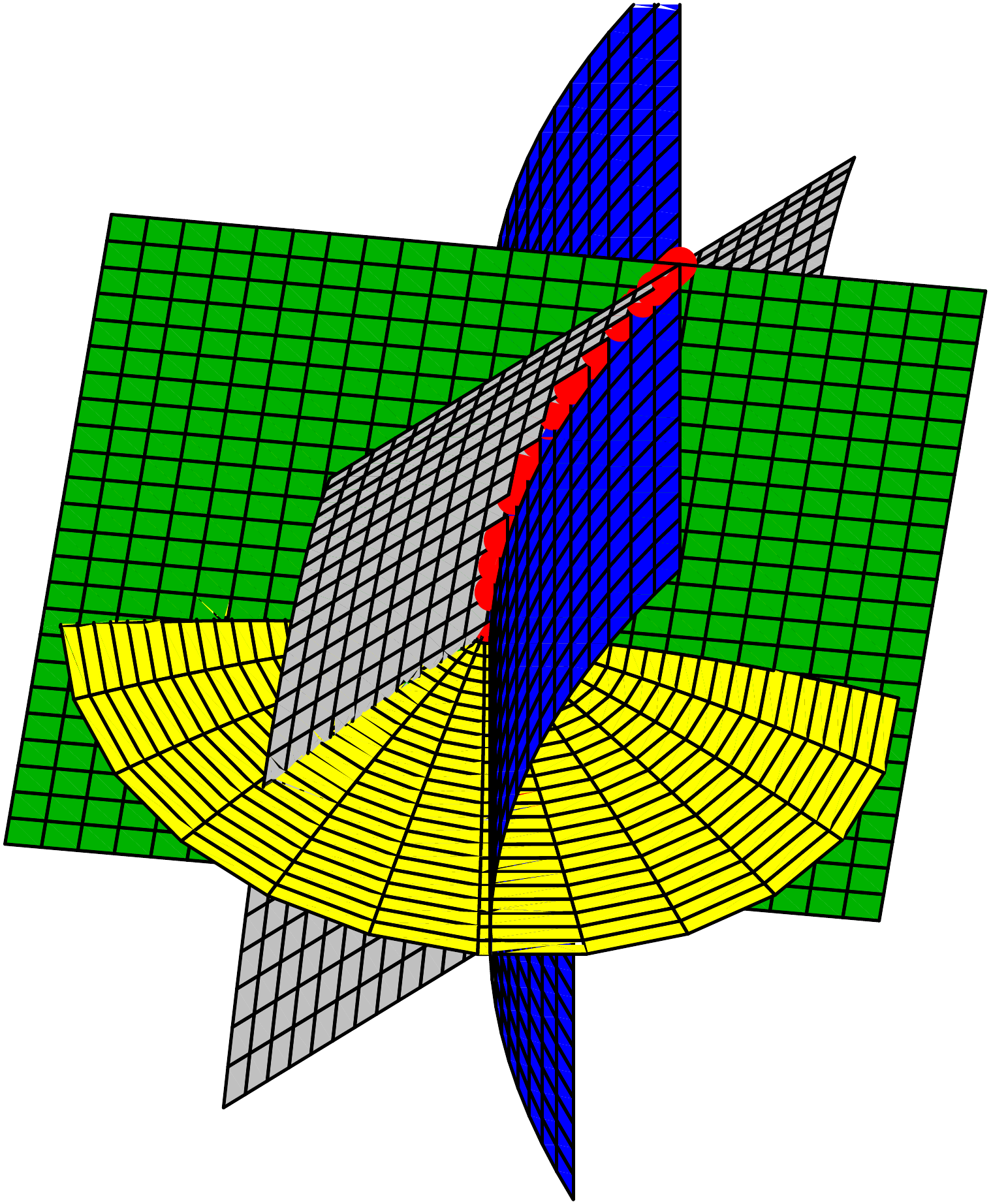}
\caption{Construction of a booklet $3$-webs from Frobenius $3$-folds.}
\label{booklet}
\end{figure}

The notion of Frobenius manifold, introduced by B.~Dubrovin (see~\cite{Df} and~\cite{Mf}), is a geometric translation of the theory of WDVV-equations that arise originally in the physical context of two-dimensional topological f\/ield theory.
\begin{definition}\label{defFrob}
A Frobenius manifold is a complex analytic manifold $M$ equipped with the following
analytic objects:
\begin{enumerate}\itemsep=0pt
\item[1)] a commutative and associative multiplication on $T_pM$,
\item[2)] an invariant non-degenerate f\/lat inner product: $\langle u \cdot v, w \rangle = \langle u, v \cdot w \rangle $,
\item[3)] a constant unity vector f\/ield $e$: $\nabla e=0$, $e \cdot v=v$ $\forall\, v\in TM$,
\item[4)] a linear Euler vector f\/ield $E$: $\nabla (\nabla E)=0$,
\end{enumerate}
satisfying the following conditions:
\begin{itemize}\itemsep=0pt
\item the f\/low of $E$ re-scales the
multiplication and the inner product,
\item 4-tensor $(\nabla_z c)(u,v,w)$ is symmetric in
$u$, $v$, $w$, $z$, where
\[
c(u,v,w):=\langle u\cdot v,w\rangle.
\]
\end{itemize}
In this def\/inition, the symbol $\nabla$ stands for the Levi-Civita connection of the inner product $\langle\ ,\ \rangle$.
\end{definition}
The following geometric construction of $3$-web via Frobenius $3$-fold was proposed in~\cite{Aw}.
Consider a {\it massive} Frobenius $3$-fold $M$, which means that the algebra $T_pM$ is semi-simple for each $p\in U$ for some open set $U\subset M$. Then $T_pM$ is the direct
product of one-dimensional algebras spanned by idempotents
\begin{equation}\label{mult}
T_pM=\mathbb C \{e_1\}\otimes \mathbb C \{e_2\}\otimes \mathbb C\{e_3\}, \qquad e_i\cdot e_j=\delta_{ij}e_i.
\end{equation}
In this setting, the unity vector f\/ield is
\[
e=e_1+e_2+e_3.
\]
Let $S$ be a surface transverse to the unit vector f\/ield $e$, then 3
planes spanned by $\{e,e_i\}$ cut 3 directions on $T_pS$ (see Fig.~\ref{booklet} on the left).
The integral curves of these direction f\/ields build a {\it flat} (or {\it hexagonal}) {\it $3$-web}, i.e., at the points with pairwise distinct web directions, this web is locally biholomorphic to 3 families of parallel lines in the plane.
Due to the existence of {\it canonical Dubrovin} coordinates (see~\cite{Di}), the distributions $\{e,e_i\}$ are integrable, integral surfaces of each of the distributions being formed by the integral curves of the unity vector f\/ield~$e$. Thus for each point in $M$ there are 3 integral surfaces intersecting along such a curve. These surfaces cut $S$ along the constructed web. This justif\/ies the following def\/inition (see also Fig.~\ref{booklet} on the right).

\begin{definition}
The constructed web is called a booklet $3$-web.
\end{definition}
The web directions of a booklet $3$-web are well-def\/ined everywhere (see~\cite{Aw}), maybe with multiplicity.
We call a point {\it regular or non-singular} if the web directions at this point are pairwise distinct. The locus of singular points is called the discriminant curve.
By def\/inition, a hexagonal $3$-web does not have any local invariants at regular points. Its ``personality'' is encoded in the behavior at \emph{singular points}, where at least 2 web directions coincide.

A booklet $3$-web enjoys the following properties~\cite{Aw}:
\begin{itemize}\itemsep=0pt
\item it has at least one inf\/initesimal symmetry at each singular point,
\item its Chern connection form is closed, i.e., the web is f\/lat,
\item its Chern connection form is holomorphic at singular points (in a suitable normalization),
\item it is biholomorphic to the characteristic $3$-web of the corresponding solution of associativity equation,
\item its Chern connection is induced by the connection on $TM$ compatible with the algebraic structure of $M$.
\end{itemize}
Note that the existence of a symmetry at a singular point is not a trivial condition for a f\/lat $3$-web: even though a f\/lat $3$-web has $3$-dimensional symmetry algebra at regular points (see~\cite{Cg}), not all symmetries survive at singular ones (see the discussion in~\cite{Ac}). For a booklet $3$-web the symmetry is generated by the f\/low of the Euler f\/ield $E$.

The analyticity condition of the Chern connection at singular points is also rather restrictive: a f\/lat $3$-web can have closed but not holomorphic connection form. Generically it has a pole on the discriminant curve (see~\cite{Ac} for examples).

For studying our webs we use implicit ODEs; this approach has proved its ef\/f\/iciency
in both non-singular (see for instance~\cite{Ha}) and singular cases~\cite{Ac,Ai}.
Let us consider singular \mbox{$3$-webs} having the f\/irst three of the above properties,
and suppose that the web directions are well def\/ined (probably, with multiplicities).
Then these directions $[dx:dy]$ satisfy
a binary dif\/ferential equation
\begin{equation}\label{binary}
K_3(x,y)dy^3+K_2(x,y)dy^2dx+K_1(x,y)dydx^2+K_0(x,y)dx^3=0,
\end{equation}
where the coef\/f\/icients $K_i$ do not vanish simultaneously.
Dividing by $dx^3$ or $dy^3$ one can write it down as an implicit cubic ODE with at least one non-vanishing coef\/f\/icient. It turns out that such ODEs can be classif\/ied. Moreover, the normal form can be ef\/fectively computed, whenever the inf\/initesimal symmetry is known.
(See Section~\ref{sing}, where we present the classif\/ication results obtained in~\cite{Ac}.)

In this paper, we address a natural question: 
whether it is possible to recover a Frobenius $3$-fold germ starting from a singular $3$-web with the properties mentioned above.
Namely we give a construction of a Frobenius $3$-fold germ for each of the normal forms.

The key observation that permitted the classif\/ication is that the inf\/initesimal symmetry operator, vanishing at the base point $(0,0)$ of the web-germ, is a dilatation in suitable coordinates
\[
X=w_1x\partial_x+w_2y\partial_y.
\]
This gives us two weights $[w_1:w_2]$ of the Euler vector f\/ield, since, for booklet $3$-webs, the operator $X$ is the projection of the Euler vector f\/ield along the unity $e$ to the surface $S$ (see~\cite{Aw}).

Consider the map $\pi_e: M\to S$ that assigns to each $p\in M$ the intersection
of its trajectory under the f\/low of $e$ with $S$.
Then the vector f\/ields $e_i$ and $v_i:=d\pi_e(e_i(p))$ are $\pi$-related: $d\pi_e(e_i(p))=v_i(\pi_e(p))$.
Therefore $[v_i,v_j]=0$.

If we have a f\/lat $3$-web, then the natural candidates for $v_i$ are the commuting web direction vector f\/ields. We recover idempotents in the form $e_i=v_i+\alpha_i\partial_t$, where $t$ is a third local coordinate on a Frobenius $3$-fold germ to be constructed.
The functions $\alpha_i$ are subjected to the following conditions:
\begin{enumerate}\itemsep=0pt
\item[1)] $\sum e_i=\partial_t,$
\item[2)] $e_i$ commute: $[e_i,e_j]=0$,
\item[3)] there is a non-degenerate f\/lat metric $g$ with $g(e_i,e_j)=0$, $i\ne j$,
\item[4)] the f\/low of $X+w_3t\partial_t$ re-scales the metric for a suitable weight $w_3$.
\end{enumerate}
(By metric we understand a non-degenerate inner product.)
It turns out that these conditions can be satisf\/ied, the analyticity of the Chern connection form being crucial for the existence of the metric $g$.
We def\/ine the multiplication by formula~(\ref{mult}). The obtained structure is that of Frobenius $3$-fold since the above conditions on $\alpha_i$ are equivalent to the associativity equation with dilatation symmetry.

Thus weighted homogeneous solutions to the associativity equation, i.e.\ invariant with respect to a dilatation symmetry, can be interpreted geometrically also in terms of the web theory. Namely they describe 4-webs $\mathfrak{W}_4$ of curves in $(\mathbb C^3,0)$ with the following properties:
\begin{itemize}\itemsep=0pt
\item $\mathfrak{W}_4$ is f\/lat,
\item one of the web foliations def\/ines a non-vanishing vector f\/ield $e$,
\item $\mathfrak{W}_4$ has a 2-dimensional symmetry algebra $\{e,E\}$, where
\[
[e,E]\wedge e=0\qquad\mbox{and}\qquad e\wedge E\not\equiv0,
\]
\item the Chern connection form of the symmetry reduction of $\mathfrak{W}_4$ by $e$ remains holomorphic at singular points.
\end{itemize}
By the symmetry reduction we understand the construction similar to the construction of booklet $3$-web via $\mathfrak{W}_4$ formed by the trajectories of $e_i$.

The webs $\mathfrak{W}_4$ are recovered with some parametric arbitrariness for particular normal forms. This indicates the need for a better understanding of the webs $\mathfrak{W}_4$.

\section{Singularities of booklet webs}\label{sing}

In this section we present the classif\/ication obtained in~\cite{Ac}. To distinguish between the normal forms, we need invariants, i.e.\ objects, preserved by local biholomorphisms.

Obviously, root multiplicity and projectivised weights $[w_1:w_2]$ are invariant.
There is a~subtler invariant. Consider the cross-ratio of the three web directions
and the direction def\/ined by the inf\/initesimal symmetry.
This function is well-def\/ined on the complement of the discriminant curve
and is constant along the trajectories of the symmetry f\/low.
Thus it is a function of a~f\/irst integral of the symmetry operator.
The limit value of this cross-ratio is our third invariant.
As the cross-ratio is dependent on the order of its arguments,
we use the fol\-lo\-wing symmetrized form: multiply cubic form~(\ref{binary})
with a 1-form vanishing on the trajectories of the symmetry group
(for normal forms it is $w_xxdy-w_yydx$),
write the resulting quartic form $a_4dy^4+4a_3dy^3dx+6a_2dy^2dx^2+4a_1dydx^3+a_0dx^4$,
compute $i:=a_0a_4-4a_1a_3+3a_2^2$ and $j:=a_4a_2a_0+2a_1a_2a_3-a_2^3-a_4a_1^2-a_0a_3^2$.
Then the invariant is $[i^3:j^2]$.
The polynomials $i$, $j$ are well-known in the classical invariant theory,
being invariants of the weights 4 and 6 respectively.

\begin{theorem}[\cite{Ac}]\label{homoexact}
Suppose ODE \eqref{binary} admits an infinitesimal symmetry $X$ vanishing at the point $(0,0)$
on the discriminant curve and the germ of the Chern connection form is exact $\gamma=d(f)$,
where $f$ is some function
germ.
Then the equation and the symmetry are equivalent to one of the
following normal forms:
\begin{alignat*}{4}
&1)\quad&&y^{m_0}p^3-p=0, \quad&& X=(2+m_0)x\partial_x+2y\partial_y,&\\
&2)\quad&&p^3+2xp+y=0, \quad&& X=2x\partial_x+3y\partial_y,&\\
&3)\quad&&\left(p-\frac23x \right)\left(p^2+\frac23xp+y-\frac{2}{9}x^2\right)=0, \qquad&& X=x\partial_x+2y\partial_y,&\\
&4)\quad&&p^3+4x\left(y-\frac{4}{9}x^3\right)p+y^2+\frac{64}{81}x^6-\frac{32}{9}yx^3=0, \qquad&& X=x\partial_x+3y\partial_y,&\\
&5)\quad&&p^3+xy^2p+\frac{2}{\sqrt{27}}\frac{x^{\frac32}y^3}{\tan\left(\frac{4}{\sqrt{3}}x^{\frac32}\right)}=0, \qquad&& X=y\partial_y, &\\
&6)\quad&&p^3+y^2p=\frac{2}{\sqrt{27}}y^3\tan\big(2\sqrt{3}x+L\big), \qquad&& X=y\partial_y, &\\
&7)\quad&&p^3+y^{3+m_0}p+y^{\frac{9+3m_0}{2}}U\left(\left[(m_0+1)\right]xy^{\frac{1+m_0}{2}}\right)=0, \qquad&& X=(1+m_0)x\partial_x-2y\partial_y,&\\
&8)\quad&&p^3+xy^{3+m_0}p-\frac{x^{\frac32}y^{\frac{9+3m_0}{2}}}{V\left(\left[\frac23(m_0+1)\right]x^{\frac32}y^{\frac{1+m_0}{2}}\right)}=0, \qquad&& X=(1+m_0)x\partial_x-3y\partial_y.&
\end{alignat*}
The function $U(T)$ is defined, with $L=-\frac{2(m_0+3)}{(m_0+1)}$ and suitable constants $C_1$, $C_2$, by the relations
\begin{alignat*}{3}
& \frac{T}{3\sqrt3L}=\frac{f'\left(-\arctan\left(\frac{3\sqrt3}{2}U\right)\right)}{f\left(-\arctan\left(\frac{3\sqrt3}{2}U\right)\right)},\qquad && f(z)=\cos^{-\mu} (z)\big[C_1P^{\mu}_{\nu}(\sin z)+ C_2Q^{\mu}_{\nu}(\sin z)\big], & \\
& f\left(-\arctan\left(\frac{3\sqrt3}{2}U(0)\right)\right)\ne 0, \qquad && f'\left(-\arctan\left(\frac{3\sqrt3}{2}U(0)\right)\right)=0. &
\end{alignat*}
The initial value of $U$ vanishes $U(0)=0$ if the number $m_0$ is even.
If $U(0)\ne 0$ one can choose $0\leq\arg(U(0)) <\pi$.

The function $V(T)$ is defined, with $L=-\frac{5m_0+17}{3(m_0+1)}$, by the relations
\[
\frac{1}{3\sqrt3 L}T=\frac{f'\left(-\arctan\left(\frac{2}{3\sqrt3}V\right)\right)}{f\left(-\arctan\left(\frac{2}{3\sqrt3}V\right)\right)},\qquad f(z)=\sin^{\mu} (z)P^{\mu}_{\nu}(\cos z).
\]
In the above formulas, $P^{\mu}_{\nu}(z)$, $Q^{\mu}_{\nu}(z)$ are Legendre's functions
for $\mu=\frac{1}{2}\big(1-\frac{1}{3L}\big)$, $\nu=\frac{1}{2}\big(\frac{1}{L}-1\big)$, $m_0$ is non-negative integer,
for the form $6)$ with $L\ne 0$ one can choose $0\leq\arg(L) <\pi$.

The weights $ [w_1:w_2 ]$, the root multiplicity and the invariant $\left[i^3:j^2\right]$ uniquely determine the normal form.
\end{theorem}
\begin{remark} The invariant $\left[i^3:j^2\right]$ assumes the value $\left[0:1\right]$ for the form $5)$,
the value $\left[1:\frac{\tan^2(L)}{-27}\right]$ for the form $6)$,
and the values $\left[1:\frac{U^2(0)}{-4}\right]$, $\left[0:1\right]$ for the forms $7)$, $8)$ respectively.
\end{remark}

\section{Metric in f\/lat coordinates}
In this section we present forms of invariant metrics for our Frobenius $3$-folds in f\/lat coordinates.
For the case when all the weights of the Euler vector f\/ields are distinct,
they were obtained for arbitrary~$n=3$ in~\cite{Df}.
Here we consider a bit weaker hypothesis that the weights of~$x$ and~$y$ are not equal.
We are looking for a constant matrix
\[g=
\begin{pmatrix}
g_{tt}&g_{tx}&g_{ty}\\
g_{tx}&g_{xx}&g_{xy}\\
g_{ty}&g_{xy}&g_{yy}
\end{pmatrix},
\]
where $g_{tt}=\langle e$, $e\rangle=\langle\partial_t,\partial_t\rangle$,
$g_{tx}=\langle\partial_t,\partial_x\rangle $, \dots, $g_{yy}=\langle\partial_y,\partial_y\rangle $, satisfying
$L_E(g)={\rm const}\cdot g$. (The local f\/low $\exp(Ea)$ re-scales the metric.)
One has $L_E(g)_{tt}=-2w_tg_{tt}$, $L_E(g)_{tx}=(-w_t-w_x)g_{tx}$, \dots, $L_E(g)_{yy}=-2w_yg_{yy}$.
Therefore the non-vanishing entries in $g$ must be of the same weight.
Observe also that for each $3$-web in the classif\/ication list the inequality $w_x\ne w_y$ holds true.

$\bullet$ {\it Metric with $\langle e,e\rangle=0$.}
As $w_{tx}\ne w_{ty}$ and the metric $g$ is non-degenerate exactly one of the the entries $g_{tx}$, $g_{ty}$
is non-vanishing. Suppose that $g_{ty}=1$. Now $g_{xx}\ne 0$ since $g$ is non-degenerate.
This implies $g_{xy}=0$, $g_{yy}=0$ and $2w_x=w_t+w_y$.
If $g_{ty}=0$ then, similarly, one normalizes $g_{tx}=1$ and obtains $g_{yy}\ne 0$, $g_{xx}=g_{xy}=0$.

Suppose $g_{xx}=\delta\ne 0$ then $\langle e_1,e_2\rangle =\langle e_1,e_3\rangle =0$ implies that the weights of $p_i$ are zero. This is not possible for our singular webs. Therefore $g_{xx}$ must vanish and the form of the matrix $g$ is
\begin{equation}\label{metric0}
g=
\begin{pmatrix}
0 & 1 & 0\\
1 & 0 & 0\\
0 & 0 & \delta
\end{pmatrix},
\end{equation}
where $\delta \ne 0$.

$\bullet$ {\it Metric with $\langle e,e\rangle \ne 0$.}
At least one of the entries $g_{tx}$, $g_{ty}$ vanishes. Let $g_{ty}=0$.
Now exactly one of the the entries $g_{xy}$, $g_{yy}$ is non-vanishing.
The entry $g_{xy}$ can not vanish; otherwise $g_{yy}\ne 0$, $g_{xx}=0$, $w_t=w_y$
and because of the non-degeneracy of $g$ we have $g_{tx}\ne 0$,
which implies $w_t=w_x$ and, f\/inally, one obtains a contradiction $w_x=w_y$.
Thus we have $g_{xy}\ne 0$, $g_{yy}=0$, which implies $g_{xx}=0$.
If $g_{tx}\ne 0$ then $w_t=w_x$. And again we get the contradiction $w_x=w_y$.
Hence one can suppose that the form of the metric is
\[g=
\begin{pmatrix}
\delta & 0 & 0\\
0 & 0 & 1 \\
0 & 1 & 0
\end{pmatrix}.
\]

\section{Idempotents}
In this section we determine idempotents $e_i$ using the following 2 facts:
\begin{enumerate}\itemsep=0pt
\item[1)] the metric, being invariant, is diagonal in the basis $\{e_1,e_2,e_3\}$,
\item[2)] the idempotents commute.
\end{enumerate}
Suppose that the following implicit cubic ODE def\/ines our $3$-web:
\begin{equation}\label{flat_cub}
p^3+S(x,y)p^2+A(x,y)p+B(x,y)=0.
\end{equation}
Let $p_1$, $p_2$, $p_3$ be the roots of~(\ref{flat_cub}) at a point
$(x,y)$ outside the discriminant curve. The following 1-forms vanish on the
solutions
\begin{gather}
\sigma_1=(p_2-p_3)(dy-p_1dx),\nonumber\\
\sigma_2=(p_3-p_1)(dy-p_2dx),\nonumber\\
\sigma_3=(p_1-p_2)(dy-p_3dx)\label{normalized forms}
\end{gather}
and satisfy
\[
\sigma_1+\sigma_2 + \sigma_3=0.
\]
Let us introduce an ``area'' form by
\[
\Omega =\sigma_1\wedge
\sigma_2=\sigma_2\wedge \sigma_3=\sigma_3\wedge
\sigma_1=(p_1-p_2)(p_2-p_3)(p_3-p_1)dy\wedge dx.
\]
 The Chern
connection form is def\/ined as (see~\cite{Be})
\[
\gamma:
=h_2\sigma_1-h_1\sigma_2=h_3\sigma_2-h_2\sigma_3=h_1\sigma_3-h_3\sigma_1,
\]
where $h_i$ verify the relations
\[
d\sigma_i=h_i\Omega.
\]
The web is f\/lat if\/f the connection form
is closed: $d(\gamma)=0$. This implies $d\sigma_i=\gamma \wedge \sigma_i$ and gives an integrating factor $k$ of the forms $\sigma_i$ as a solution to the following equation
\begin{gather*}
dk=-\gamma k.
\end{gather*}
Computing the Chern connection form in terms of roots $p_i$ and using the Viete formulas one gets
\[
\gamma =\frac{\gamma_1dx+\gamma_2dy}{-D},
\]
where
\begin{gather*}
\gamma_1=\big(4BS^2-3AB-SA^2\big)S_x+B(SA-9B)S_y+(2A^2-6SB)A_x\nonumber\\
\hphantom{\gamma_1=}{}
+2B(3A-S^2)A_y+(9B-SA)B_x+\big(AS^2-4A^2+3SB\big)B_y,\nonumber\\
\gamma_2=\big(6SB-2A^2\big)S_x-2B\big(3A-S^2\big)S_y+(SA-9B)A_x\nonumber\\
\hphantom{\gamma_2=}{}
 + \big(4A^2-AS^2-3SB\big)A_y+\big(6A-2S^2\big)B_x+\big(2S^3+18B-8SA\big)B_y,
\end{gather*}
and $D$ is the discriminant of~(\ref{flat_cub}) with respect to $p$
\[
D=18SAB+S^2A^2-4A^3-27B^2-4BS^3.
\]
We also need to know how the connection form in the normalization~(\ref{normalized forms})
is transformed when we change the local coordinates
\[
\bar{y}=f(x,y), \qquad \bar{x}=g(x,y).
\]
One f\/inds without dif\/f\/iculty that the forms $\bar{\sigma_i}$
in the new coordinates are related to the forms in old ones by
\begin{gather*}
\bar{\sigma_i}=\frac{(f_yg_x-f_xg_y)^2\sigma_i}{g_x^3-Sg_x^2g_y+Ag_xg_y^2-Bg_y^3}.
\end{gather*}
Therefore (see~\cite{Be})
\begin{equation}\label{connection transforms}
\bar{\gamma}=\gamma+d\ln\left(\frac{(f_yg_x-f_xg_y)^2}{g_x^3-Sg_x^2g_y+Ag_xg_y^2-Bg_y^3}\right).
\end{equation}
Let $v_i=\xi_i \partial_x+\eta_i\partial_y$ be commuting multi-valued vector f\/ields,
whose trajectories are the web leaves. They may have singularities at singular points:
with $K$ satisfying $\frac{dK}{K}=\gamma $ we have
\begin{gather*}
v_1=\frac{K(\partial_x +p_1 \partial_y)}{(p_3-p_1)(p_1-p_2)},\qquad
v_2=\frac{K(\partial_x +p_2 \partial_y)}{(p_1-p_2)(p_2-p_3)},\qquad
v_3=\frac{K(\partial_x +p_3 \partial_y)}{(p_2-p_3)(p_3-p_1)}.
\end{gather*}
Here $x$, $y$ are restrictions on $S$ of f\/lat coordinates in $M$, in which the Euler vector f\/ield has the form $E=w_tt\partial_t+w_xx\partial_x+w_yy\partial_y$. A subtle point is that f\/lat coordinates are not necessarily the coordinates used for the normal forms: they are related to them by some coordinate transform preserving $X$.

The idempotents are
\begin{equation}\label{idempotents}
e_1=\alpha \partial_t+v_1, \qquad e_2=\beta \partial_t+v_2,\qquad e_1=(1-\alpha -\beta) \partial_t+v_3,
\end{equation}
where the components $\alpha$, $\beta$ are to be def\/ined by the orthogonality condition
\begin{equation}\label{orthogonal}
\langle e_1,e_2\rangle =\langle e_2,e_3\rangle =\langle e_3,e_1\rangle =0.
\end{equation}

$\bullet$ {\it Metric with $\langle e,e\rangle =0$.}
Equations~(\ref{orthogonal}) give
\begin{equation}\label{dfor0}
\delta=-1/K.
\end{equation} Adjusting the integration constant so that $K=-1$ we have $\delta =1$ and
\begin{equation}\label{alphabeta0}
\alpha =\frac{p_2p_3-p_1(p_2+p_3)}{2(p_1-p_2)(p_1-p_3)},\qquad \beta =\frac{p_1p_3-p_2(p_1+p_3)}{2(p_2-p_1)(p_2-p_3)}.
\end{equation}
Since $dK=0$, we see that $\gamma_1=\gamma_2=0$.
The equation $[e_1,e_2]=0$ reads as $v_1(\beta)=v_2(\alpha)$ and obviously implies $[e_1,e_3]=[e_2,e_3]=0$.
The above 3 equations $\gamma_1=\gamma_2=v_1(\beta)-v_2(\alpha)=0$
can be resolved with respect to $S_x$, $A_x$, $B_x$
to give the following system of hydrodynamic type
\begin{gather}\label{WDVV0}
S_x= \frac{1}{2}A_y,\qquad
A_x= 2By,\qquad
B_x= SB_y+BS_y-\frac{1}{2}AA_y.
\end{gather}
\begin{lemma}\label{web_ass_0}
Suppose that the solution web of ODE~\eqref{flat_cub} satisfies the following conditions:
\begin{enumerate}\itemsep=0pt
\item[$1)$] it admits an infinitesimal symmetry $X=w_xx\partial_x+w_yy\partial_y$ with $w_x\ne w_y$,
\item[$2)$] the coefficients $S$, $A$, $B$ satisfy the system of PDEs~\eqref{WDVV0}.
\end{enumerate}
Then there is a Frobenius $3$-fold germ with $\langle e,e\rangle =0$ whose booklet web is the solution web of~\eqref{flat_cub}.
\end{lemma}
\begin{proof}
Equations~(\ref{WDVV0}) ensures the local existence of a function $f(x,y)$ satisfying
$f_{yyy}=S$, $f_{yyx}=\frac{1}{2}A$, $f_{yxx}=B$, $f_{xxx}=BS-\frac{1}{4}A^2$.
Thus we obtain the associativity equation
\begin{equation}\label{assPoly}
f_{xxx}=f_{yyy}f_{yxx}-f^2_{yyx}.
\end{equation}
(See~\cite{MFa} where the equivalence of~(\ref{WDVV0}) and~(\ref{assPoly}) was established and~\cite{Al} where equations of these types were studied.)
Since the web admits an inf\/initesimal symmetry $X$, the function~$f$ is weighted homogeneous and def\/ines a germ of Frobenius $3$-fold with $\langle e,e\rangle =0$ (see~\cite{Df}). Equation~(\ref{flat_cub}) def\/ines the corresponding characteristic web of the above associativity equation. Hence the web is the booklet web of the constructed Frobenius $3$-fold germ (see~\cite{Aw}).
\end{proof}
The value of the lemma is that it makes unnecessary the checking the potentiality condition for the 4-tensor $(\nabla_z c)(u,v,w)$ in Def\/inition~\ref{defFrob}.
\begin{theorem}\label{web_Frobenius_0}
Suppose that the solution web of ODE~\eqref{flat_cub} satisfies the following conditions:
\begin{enumerate}\itemsep=0pt
\item[$1)$] it admits an infinitesimal symmetry $X=w_xx\partial_x+w_yy\partial_y$ with $w_x\ne w_y$,
\item[$2)$] its Chern connection form vanishes identically $\gamma=0$,
\item[$3)$] the vector fields $e_1$, $e_2$ commute, where $\alpha$, $\beta$ are as in~\eqref{alphabeta0}.
\end{enumerate}
Then there is a Frobenius $3$-fold germ with $\langle e,e\rangle=0$
whose booklet web is the solution web of~\eqref{flat_cub}.
\end{theorem}
\begin{proof}
Def\/ine a metric by~(\ref{metric0}) with $\delta=1$,
the multiplication as the direct sum of one-dimensional algebras with unities~(\ref{idempotents}),
and the Euler vector f\/ield by $E=w_tt\partial_t+X$ with $w_t=2w_y-w_x$.
Then Lemma~\ref{web_ass_0} ensures that the def\/ined structure is that of Frobenius \mbox{$3$-fold}.
\end{proof}

$\bullet$ {\it Metric with $\langle e,e\rangle \ne 0$.}
Equations~(\ref{orthogonal}) give
\begin{equation}\label{alphabeta1}
\alpha =\frac{(p_1+p_2)(p_1+p_3)}{(p_1-p_2)(p_1-p_3)},\qquad \beta = \frac{(p_2+p_1)(p_2+p_3)}{(p_2-p_1)(p_2-p_3)},
\end{equation}
and
\[
\delta=\frac{-K^2}{(p_1+p_2)(p_2+p_3)(p_3+p_1)}.
\]
 Vieta formulas result in
\begin{equation}\label{dfor1}
\delta=\frac{K^2}{SA-B}.
\end{equation}
Adjusting the integration constant so that $K=-1$ we have $\delta =1$, provided the condition $2\frac{dK}{K}=\frac{d(SA-B)}{SA-B}$ is satisf\/ied. Substituting $\frac{dK}{K}=\gamma$ we obtain
\begin{equation}\label{connection condition for 1}
2\gamma=\frac{d(SA-B)}{SA-B}.
\end{equation}
Together with the commutativity condition $v_1(\beta)=v_2(\alpha)$, where $\alpha$, $\beta$ are given now by~(\ref{alphabeta1}), equation~(\ref{connection condition for 1}) gives the following system of hydrodynamic type
\begin{gather}
S_x= \frac{(A^2-AS^2+2SB)S_y+(S^3-AS+2B)A_y-(A+S^2)B_y}{2AS-2B},\nonumber\\
A_x= \frac{-AS_y+SA_y+B_y}{2},\nonumber\\
B_x= \frac{(A^3+4B^2-3ASB)S_y+(2AB-SA^2+BS^2)A_y+(2AS^2-3SB-A^2)B_y}{2AS-2B}.\label{WDVV1}
\end{gather}
\begin{lemma}\label{web_ass_1}
Suppose that the solution web of ODE~\eqref{flat_cub} satisfies the following conditions:
\begin{enumerate}\itemsep=0pt
\item[$1)$] it admits an infinitesimal symmetry $X=w_xx\partial_x+w_yy\partial_y$ with $w_x\ne w_y$,
\item[$2)$] the coefficients $S$, $A$, $B$ of~\eqref{flat_cub} satisfy the system of PDEs~\eqref{WDVV1}.
\end{enumerate}
Then there is a Frobenius $3$-fold germ with $\langle e,e\rangle =1$
whose booklet web is the solution web of~\eqref{flat_cub}.
\end{lemma}
\begin{proof} It repeats that of Lemma~\ref{web_ass_0}.
System~\eqref{WDVV1} again is reducible (see~\cite{Al})
and therefore is equivalent to the local existence of a function $f(x,y)$
satisfying $\frac{f_{yyx}}{f_{yyy}}=S$, $-\frac{f_{yxx}}{f_{yyy}}=A$, $-\frac{f_{xxx}}{f_{yyy}}=B$,
and the associativity equation
\begin{equation}\label{assDet}
f_{xxx}f_{yyy}-f_{xxy}f_{xyy}=1,
\end{equation}
 which corresponds to the case $\langle e,e\rangle =1$ (see~\cite{Df}).
\end{proof}
\begin{theorem}\label{web_Frobenius_1}
Suppose that the solution web of ODE~\eqref{flat_cub} satisfies the following conditions:
\begin{enumerate}\itemsep=0pt
\item[$1)$] it admits an infinitesimal symmetry $X=w_xx\partial_x+w_yy\partial_y$ with $w_x\ne w_y$,
\item[$2)$] its Chern connection form verifies~\eqref{connection condition for 1},
\item[$3)$] the vector fields $e_1$, $e_2$ commute, where $\alpha$, $\beta$ are as in~\eqref{alphabeta1}.
\end{enumerate}
Then there is a Frobenius $3$-fold germ with $\langle e,e\rangle =1$ whose booklet web is the solution web of~\eqref{flat_cub}.
\end{theorem}
\begin{proof} The proof is similar to that of Theorem~\ref{web_Frobenius_0}.
The only dif\/ference is that the weight $w_t$ now is $w_t=(w_x+w_y)/2$.
\end{proof}
\begin{remark}
Equations~(\ref{dfor0}) and~(\ref{dfor1}) imply that the Chern connection form of the booklet web remains holomorphic at singular points: this is straightforward for equation~(\ref{dfor0}), for~(\ref{dfor1}) one has to take into account the associativity equation~(\ref{assDet}) and, if equation~(\ref{binary}) is not monic, the formula~(\ref{connection transforms}). Thus, a geometrical interpretation of the analyticity of the connection is the existence of non-degenerate invariant f\/lat metric.
\end{remark}

\section{Webs with one elliptic symmetry operator}
We call a symmetry operator $X=w_xx\partial_x+w_yy\partial_y$ elliptic if the weight ratio $w_x/w_y$ is positive. In this section we study the Frobenius $3$-fold germs corresponding to the normal forms 2), 3) and 4) in Theorem~\ref{homoexact}. The connection form $\gamma$ vanishes identically for each of these equations.
The coordinate transformations, preserving the symmetry generator $X$, are non-trivial only for the forms 3) and 4). Up to an unessential scaling, they read as
\begin{equation}\label{transform3}
\bar{x}=x, \qquad \bar{y}=y+rx^2,\qquad r=\operatorname{const},
\end{equation}
and
\begin{equation}\label{transform4}
\bar{x}=x, \qquad \bar{y}=y+rx^3,\qquad r=\operatorname{const},
\end{equation}
respectively.
\begin{proposition}
The booklet $3$-web of a massive Frobenius $3$-fold with $\langle e,e\rangle =1$ cannot have singularities of types
$2)$, $3)$, $4)$ in Theorem~{\rm\ref{homoexact}}.
\end{proposition}
\begin{proof}
The coef\/f\/icients $\bar{S}$, $\bar{A}$, $\bar{B}$ of the transformed ODEs vanish at singular points.
The condition $\gamma \equiv 0$ is also preserved by transformations~(\ref{transform3}) and~(\ref{transform4}).
This yields a contradiction: $\frac{1}{\delta}=\frac{\bar{S}\bar{A}-\bar{B}}{K^2}=0$.
\end{proof}
\begin{definition}[\cite{Df}]
Two Frobenius manifolds $M$ and $\tilde{M}$ are equivalent if there exists a dif\/feo\-mor\-phism
\begin{gather*}
 \varphi: \  M \to \tilde{M}
\end{gather*}
being a linear conformal transformation of the corresponding invariant metrics
\begin{gather*}
\varphi ^*\tilde{g}=c^2g
\end{gather*}
($c$~is a non-zero constant) with the dif\/ferential acting as an isomorphism on the tangent algebras
\begin{gather*}
d\varphi: \  T_pM \to T_{\varphi (p)}\tilde{M}.
\end{gather*}
\end{definition}
\begin{proposition}
Each of the $3$-web germs of type $2)$, $3)$, $4)$
in Theorem~{\rm \ref{homoexact}} is a booklet $3$-web for some Frobenius $3$-fold germ
with $\langle e,e\rangle =0$. Moreover, this germ is unique up to the equivalence of Frobenius $3$-folds.
\end{proposition}
\begin{proof} We use Lemma~\ref{web_ass_0}. The coef\/f\/icients of ODE $2)$ satisfy conditions~(\ref{WDVV0}).
Now let us apply transformation~(\ref{transform3}) to ODE $3)$ and transformation~(\ref{transform4}) to ODE $4)$.
Then the f\/irst equation of~(\ref{WDVV0}) gives $r=-\frac{1}{12}$ for $3)$ and $r=-\frac{1}{9}$ for $4)$,
the second and the third equations~(\ref{WDVV0}) being satisf\/ied for the found values of~$r$.
The constructed germ is obviously unique up to an unessential scaling, i.e., up to the equivalence.
\end{proof}
\begin{remark}
The Frobenius $3$-folds in the above Proposition are those corresponding to the polynomial solutions of WDVV equation~(\ref{assPoly}) (see~\cite{Df}), which, of course, is not surprising. The above Frobenius $3$-fold germ can be also recovered directly via Theorem~\ref{web_Frobenius_0}.
\end{remark}

\section{Web with a parabolic symmetry operator}
We call a symmetry operator $X=w_xx\partial_x+w_yy\partial_y$ parabolic if one of the weights $w_x$, $w_y$ vanishes.
In this section we study the Frobenius $3$-fold germs corresponding to the normal forms~5) and~6)
of Theorem~\ref{homoexact}.
Coordinate transformations, preserving the symmetry generator $X$, are of the form
\begin{gather*}
\bar{y}=yF(x), \qquad \bar{x}=G(x),
\end{gather*}
where
\begin{gather*}
F(0)\ne 0, \qquad G'(0)\ne 0, \qquad G(0)=0.
\end{gather*}

\begin{proposition}
The booklet $3$-web of a massive Frobenius $3$-fold with $\langle e,e\rangle =1$
cannot have singularities of type $5)$ or {\rm 6)} of Theorem~{\rm\ref{homoexact}}.
\end{proposition}
\begin{proof}
The orthogonality condition~(\ref{orthogonal}) yields the contradiction
\begin{gather*}
\frac{1}{\delta}=\left.-\frac{BF^3+2yAF^2 F'+8y^3(F')^3}{K^2F^4G'}\right|_{(0,0)}=0.\tag*{\qed}
\end{gather*}
\renewcommand{\qed}{}
\end{proof}

\begin{proposition}
 There is a one-parameter family of inequivalent Frobenius $3$-fold germs with $\langle e,e\rangle =0$,
 whose booklet web germs are diffeomorphic to $3$-web germ $5)$ or $6)$ of Theorem~{\rm\ref{homoexact}}.
\end{proposition}
\begin{proof} Consider the following cubic ODE
\[
p^3+ys(x)p^2+y^2a(x)p+y^3b(x)=0.
\]
Equations~(\ref{WDVV0}) read as
\[
s'=a,\qquad a'=6b, \qquad b'=4ab-a^2.
\]
The invariant $[i^3:j^2]$ is equal to $\big[108(s^2-3a)^3:(9as-2s^3-27b)^2\big]$.
Adjusting the initial conditions $s(0)$, $a(0)$, $b(0)$ to verify
\[
\left.\frac{(9as-2s^3-27b)^2}{(3a-s^2)^3}\right|_{x=0}=4\tan^2(L)
\]
 and applying Lemma~\ref{web_ass_0} one gets 2-parameter family Frobenius $3$-fold germs with the booklet web 6).
For the form $5)$ one chooses $\left.s^2-3a\right|_{x=0}=0$ and $\left.9as-2s^3-27b\right|_{x=0}\ne 0$. Now the symmetry $y\partial_y$ reduces the number of parameters in the family of inequivalent Frobenius $3$-fold germs to one.
\end{proof}
\begin{remark} Frobenius $3$-fold germ can be also constructed directly via Theorem~\ref{web_Frobenius_0}. Consider, for instance, the case 6).
The connection form here is
$\gamma=\frac{2}{\sqrt3}\tan(2\sqrt3 x+L)dx,$ therefore $K=\frac{1}{\sqrt[3]{\cos(2\sqrt 3 x+L)}}$.
 The orthogonality conditions~(\ref{orthogonal}) give
 \begin{equation}\label{dG1F}
 \delta= -\frac{G'}{KF^2},
 \end{equation}
 and the following expressions for the idempotent components
 \[
 \alpha=\frac{(p_2p_3-p_1p_2-p_1p_3)F^2-2yp_1FF'-y^2(F')^2}{2F^2(p_1-p_2)(p_1-p_3)},
 \]
 \[
 \beta=\frac{(p_1p_3-p_2p_1-p_2p_3)F^2-2yp_2FF'-y^2(F')^2}{2F^2(p_2-p_1)(p_2-p_3)}.
 \]
 Now the commutativity condition
$[e_1,e_2]=0$ amounts to the following equation for $F$
 \begin{equation}\label{Fparabolic}
 3FF''-6(F')^2-2\sqrt3\tan\big(2\sqrt3x+L\big)FF'+F^2=0.
 \end{equation}
 Due to Theorem~\ref{web_Frobenius_0}, each solution $G$, $F$ of equations~(\ref{dG1F}) and~(\ref{Fparabolic})
 def\/ines a Frobenius $3$-fold germ.
\end{remark}
\begin{remark}
Substitution $\frac{F'}{F}=u$ reduces equation~(\ref{Fparabolic}) to some Riccati equation.
Hence~(\ref{Fparabolic}) is linearizable.
\end{remark}

\section{Web with a hyperbolic symmetry operator}
We call a symmetry operator $X=w_xx\partial_x+w_yy\partial_y$ hyperbolic if the weight ratio $w_x/w_y$ is negative. In this section we study the Frobenius $3$-fold germs corresponding to the normal forms $7)$ and $8)$ of Theorem~\ref{homoexact}.
The connection form here is holomorphic, i.e., the factor $K$ does not vanish.
Coordinate transformations, preserving the symmetry generator $X$, are of the form
\begin{equation}\label{transform6}
\bar{y}=yQ(s), \qquad \bar{x}=xR(s),
\end{equation}
 where $Q(0)\ne 0$, $R(0)\ne 0$, $s=xy^r$, $r=\frac{1+m_0}{2}$, and the functions $Q$, $R$ are even for even $m_0$.

\begin{proposition}
The booklet $3$-webs of a massive Frobenius $3$-fold with $\langle e,e\rangle =1$
cannot have singularities of types $7)$ or $8)$ of Theorem~{\rm\ref{homoexact}}.
\end{proposition}
\begin{proof} Since the Jacobian matrix of transformation~\eqref{transform6} is diagonal at $(0,0)$,
we infer that the transformed cubic ODE remains monic 
with coef\/f\/icients $\bar{S}$, $\bar{A}$, $\bar{B}$ vanishing at $(0,0)$.
Further, the factor $\bar{K}$ remains analytic and f\/inite due to~(\ref{connection transforms}).
Now the condition~(\ref{dfor1})
gives a contradiction
\[
\frac{1}{\delta}=\frac{\bar{S}\bar{A}-\bar{B}}{\bar{K}^2}=0.\tag*{\qed}
\]
\renewcommand{\qed}{}
\end{proof}

To prove that each singularity of types $7)$ and $8)$ can be realized by some booklet web,
we need some properties of the coef\/f\/icients $\bar{S}$, $\bar{A}$, $\bar{B}$.
One easily f\/inds by direct calculation the following formulae for them
\begin{equation}\label{SAB}
\bar{S}=\frac{y}{x}\sigma(s), \qquad \bar{A}=\frac{y^2}{x^2}\alpha(s), \qquad \bar{B}=\frac{y^3}{x^3}\beta(s),
\end{equation}
where
\begin{gather}
\sigma(s)=-\frac{s}{\Delta}\left[3s^2\big(1+r^2s^2+r^3Fs^3\big)Q'(R')^2+rs^3(2+3rFs)Q(R')^2\right.\nonumber\\
\left. \phantom{\sigma(s)=}{} +2s\big(3+r^2s^2\big)RQ'R'+2rs^2RQR'+3R^2Q'\right],\nonumber\\
\alpha(s)=\frac{s^2}{\Delta} \left[RQ^2+\big(3+r^2s^2\big)R(Q')^2+3s\big(1+r^2s^2+r^3Fs^3\big)(Q')^2R'\right.\nonumber\\
\left. \phantom{\alpha(s)=}{} +2rs^2(2+3rFs)QQ'R'+s(1+3rFs)Q^2R'+2rsRQQ'\right],\nonumber\\
\beta(s)=-\frac{s^3}{\Delta} \left[ rs(2+3rFs)Q(Q')^2+\big(1+r^2s^2+r^3Fs^3\big)(Q')^3+(1+3rFs)Q^2Q'+FQ^3\right],\nonumber\\
\Delta=R^3+3sR^2R'+s^2\big(3+r^2s^2\big)R(R')^2+s^3\big(1+r^2s^2+r^3Fs^3\big)(R')^3.\label{sab}
\end{gather}
\begin{proposition}
For each $3$-web germ singularity of type $7)$ with odd $m_0$ or type $8)$ of Theorem~{\rm\ref{homoexact}}
there is a one-parameter family of inequivalent Frobenius $3$-fold germs with $\langle e,e\rangle =0$, whose booklet web germs are diffeomorphic to this web. For each $3$-web germ singularity of type $7)$ with even $m_0$ there is, up to equivalence, unique Frobenius $3$-fold germ with $\langle e,e\rangle =0$.
\end{proposition}
\begin{proof} Theorem~\ref{web_Frobenius_0} is not of much use here: calculations become very involved.
We prove the existence using Lemma~\ref{web_ass_0}. Coef\/f\/icients $S$, $A$, $B$ of equation~\eqref{flat_cub} have the form~(\ref{SAB}). Then equations~(\ref{sab}) imply
\[
\sigma(s)=s\tilde{\sigma}(s), \qquad \alpha(s)=s^2\tilde{\alpha}(s), \qquad \beta(s)=s^3\tilde{\beta}(s).
\]
Now the system~(\ref{WDVV0}) gives
\begin{gather}
\tilde{\sigma}'=\frac{(k+1)\tilde{\alpha}-k(k+1)s(\tilde{\sigma}\tilde{\alpha}-3\tilde{\beta})+k^2(k+1)s^2\tilde{\sigma}\tilde{\beta}}
{1-ks\tilde{\sigma}+k^2s^2\tilde{\alpha}-k^3s^3\tilde{\beta}},\nonumber\\
\tilde{\alpha}'=\frac{6(k+1)\tilde{\beta}-2k(k+1)s(\tilde{\alpha}^2-\tilde{\sigma}\tilde{\beta})+2k^2(k+1)s^2\tilde{\alpha}\tilde{\beta}}
{1-ks\tilde{\sigma}+k^2s^2\tilde{\alpha}-k^3s^3\tilde{\beta}},\nonumber\\
\tilde{\beta}'=\frac{-(k+1)(\tilde{\alpha}^2-4\tilde{\sigma}\tilde{\beta})-2k(k+1)s\tilde{\alpha}\tilde{\beta}+3k^2(k+1)s^2\tilde{\beta}^2}
{1-ks\tilde{\sigma}+k^2s^2\tilde{\alpha}-k^3s^3\tilde{\beta}}.\label{eq_for_sab}
\end{gather}
This system has solutions for each choice of initial conditions $\tilde{\sigma}(0)$, $\tilde{\alpha}(0)$, $\tilde{\beta}(0)$.
In terms of $\tilde{\sigma}(0)$, $\tilde{\alpha}(0)$, $\tilde{\beta}(0)$ the invariant $ [i^3:j^2 ]$ reads as
\[
\Big[108\big\{\tilde{\sigma}^2(0)-3\tilde{\alpha}(0)\big\}^3:\big\{2\tilde{\sigma}^3(0)-9\tilde{\alpha}(0)\tilde{\sigma}(0)+27\tilde{\beta}(0)\big\}^2\Big].
\]
For the form $8)$ we set $\tilde{\alpha}(0)=\tilde{\sigma}^2(0)/3$,
and choose $\tilde{\sigma}(0)$, $\tilde{\beta}(0)$ to satisfy $\tilde{\sigma}^3(0)\ne 27\tilde{\beta}^2(0)$.
The symmetry generated by the Euler vector f\/ield reduces the number of free parameters by one, therefore we get one-parameter family of inequivalent $3$-fold germs.

For the form $7)$ consider f\/irst the case of odd $m_0$.
Adjusting the initial conditions so that the value of the invariant coincides with the corresponding value of the normal form, and taking into account the symmetry, we obtain again one-parameter family of inequivalent $3$-fold germs.

Now suppose that $m_0$ is even. System~(\ref{eq_for_sab})
is symmetric with respect to the following involution:
$s\to -s$, $\tilde{\sigma}\to -\tilde{\sigma}$, $\tilde{\alpha}\to \tilde{\alpha}$,
$\tilde{\beta}\to -\tilde{\beta}$.
Uniqueness of the solution of Cauchy problem with initial conditions
$\tilde{\sigma}(0)=0$, $\tilde{\alpha}(0)=a_0$, $\tilde{\beta}(0)=0$
ensures that the function $\tilde{\sigma}$ and $\tilde{\beta}$ are odd, while~$\tilde{\alpha}$ is even.
Therefore the coef\/f\/icients
\begin{gather*}
S=\frac{y}{x}s\tilde{\sigma}(s), \qquad A=\frac{y^2}{x^2}s^2\tilde{\alpha}(s), \qquad B=\frac{y^3}{x^3}s^3\tilde{\beta}(s)
\end{gather*}
of ODE~\eqref{flat_cub} are holomorphic and satisfy~(\ref{WDVV0}). Here the structure of Frobenius $3$-fold is unique up to the equivalence generated by the Euler vector f\/ield.
\end{proof}

\section{Web with 2-dimensional symmetry algebra}
In this section we consider the form 1).
Here the symmetry algebra is 2-dimensional:
note that the equation 1) is invariant with respect to the translation operator $\partial_x$.
For a Frobenius $3$-fold to be constructed we have an alternative:
either the Euler vector f\/ield vanishes at $(0,0)$ or it does not vanish.

\begin{proposition}
The booklet $3$-web germ of a massive Frobenius $3$-fold cannot have singularities of type $1)$
of Theorem~{\rm\ref{homoexact}},
if the projection of the Euler vector field to the surface $S$ along the unity vector $e$ vanishes.
\end{proposition}
\begin{proof}
The symmetry operator giving rise to the Euler vector f\/ield is
\[
(1+{m_0}/{2})x\partial_x+y\partial_y.
\]
Then the f\/lat coordinates are, up to an unessential scaling
\[
\bar{x}=x+cy^{1+{m_0}/{2}}, \qquad \bar{y}=y.
\]
The orthogonality condition~(\ref{orthogonal}) yields: $\delta=0$ for $\langle e,e\rangle=0$
and $\frac{1}{\delta}=0$ for $\langle e,e\rangle \ne 0$.
\end{proof}

Unfortunately, because of 
computational dif\/f\/iculties,
the author is not able to provide a~complete analysis for the case when the projection of $E$ to $S$ does not vanish.

\begin{proposition}
For each $3$-web germ singularity of type $1)$
there is parametric family of inequivalent Frobenius $3$-fold germs with $\langle e,e\rangle =0$, whose booklet web germs are diffeomorphic to this web.
\end{proposition}
\begin{proof} The proof is somewhat of an ``experimental'' nature: the author used the symbolic computational system Maple.

 Let us choose the symmetry operator as $X=(1+(l+1)x)\partial_x+y\partial_y$ and $1+(l+1)x$ as a~new coordinate $z$.
 Now the base point of the $3$-fold germ is $(1,0)$. We are looking for the f\/lat coordinates in the form
\[
\bar{x}=z^{\alpha}R(t), \qquad \bar{y}=z^{\beta}Q(t),\qquad t=yz^{-\frac{1}{1+m_0/2}}
\]
with analytic $R$ and $Q$ subjected to $\left.\alpha RQ'-\beta QR'\right|_0\ne 0$.

The orthogonality condition~(\ref{orthogonal}) implies that necessarily $\alpha=2\beta -1$ for $\langle e,e\rangle =0$.
The analysis of ODEs for $R$ and $S$, arising from the hypotheses of Theorem~\ref{web_Frobenius_0}, shows that they have at least local solutions at $t=0$ for any $\beta \ne \frac{1}{2}$. Therefore one can construct at least 1-parametric family of inequivalent Frobenius $3$-fold germs.
\end{proof}
\begin{remark}
The symmetry operator in the f\/lat coordinates is $(2\beta -1)z\partial_z+\beta y\partial_y$.
Observe that~$\beta$ is not necessarily rational, unlike the case of triple root.
\end{remark}

\section{Concluding remarks}

{\bf 1.} The 4-web of curves $\mathfrak{W}_4$ discussed in Introduction is formed by the trajectories of the unity vector f\/ield $e$ and of the vector f\/ields of idempotents $e_i$. The generalization to Frobenius manifolds of higher dimensions is straightforward. We hope that the study of this object, namely, $(n+1)$-web of curves $\mathfrak{W}_{n+1}$ in $n$-dimensional Frobenius manifold $M^n$ will provide a better insight into the singular set (or ``discriminant'') of the Frobenius manifold. Note that for the $3$-webs with parabolic and hyperbolic symmetries the 4-web $\mathfrak{W}_4$ are recovered with some parametric arbitrariness. We hope to interpret these parameters in terms of $\mathfrak{W}_4$.

It seems that the web of curves, with the exception of webs in the plane, was not a very popular object to study in dif\/ferential geometry (see surveys~\cite{GSr} and~\cite{AGh}).

{\bf 2.} The reduction of $\mathfrak{W}_{n+1}$ by the symmetry generated by the unity vector f\/ield $e$ can also be generalized to higher dimensions. As a result we obtain, for $n$-dimensional Frobenius manifold, a f\/lat $n$-web germ of curves in $(\mathbb C^{n-1},0)$ admitting a ``linear'' symmetry.

{\bf 3.} It is interesting that Frobenius $3$-folds with $\langle e,e\rangle \ne 0$ does not produce singular booklet $3$-webs with triple root and vanishing symmetry. The author does not know if one can get a~singular booklet $3$-webs with non-vanishing symmetry at this situation. An approach based on Lemma~\ref{web_ass_1} gives singular ODEs for the coef\/f\/icients
$S=\frac{y}{x}\tilde{\sigma}(s)$, $A=\frac{y^2}{x^2}\tilde{\alpha}(s)$, and $B=\frac{y^3}{x^3}\tilde{\beta}(s)$, if one imposes initial conditions corresponding to the desired singularity.

\subsection*{Acknowledgements}

This research was partially supported by MCT/CNPq/MEC/CAPES~-- Grant 552758/2011-6.

\pdfbookmark[1]{References}{ref}
\LastPageEnding

\end{document}